\documentclass{amsart}
\usepackage{amsfonts}
\usepackage{amsmath}
\usepackage{amssymb}
\usepackage{amstext}

\def\rat{{\mathbb Q}}
\def\integ{{\mathbb Z}}
\def\calo{{\mathcal O}}
\def\calt{{\mathcal T}}
\def\caln{{\mathcal N}}

\newcommand{\invlim}[1]{\lim_{\stackrel{\leftarrow}{#1}}}
\def\tensor{\otimes}

\def\grh{{{\rm (GRH)}}}
\def\nogrh{{{\rm (U)}}}
\def\lagod{{\rm (LO)}}
\def\bach{{\rm (BS)}}

\def\cross{\times}
\def\units{^\cross}

\def\idp{{\mathfrak p}}
\def\idq{{\mathfrak q}}

\newcommand{\artin}[2]{\left[\frac{#2}{#1}\right]}
\newcommand{\til}[1]{{\widetilde{#1}}}
\newcommand{\st}[1]{\{#1\}}
\newcommand{\abs}[1]{{\left|#1\right|}}

\DeclareMathOperator{\gsp}{GSp}
\DeclareMathOperator{\gal}{Gal}
\DeclareMathOperator{\frob}{Fr}
\DeclareMathOperator{\End}{End}
\DeclareMathOperator{\gl}{GL}
\DeclareMathOperator{\aut}{Aut}
\DeclareMathOperator{\tr}{tr}
\DeclareMathOperator{\Frac}{Frac}

\def\iso{\cong}
\def\defeq{:=}
\def\inject{\hookrightarrow}
\def\ra{\rightarrow}

\def\norm{{\bf N}}
\newcommand{\disc}[1]{\Delta_{{#1}}}

\def\jcite{\cite}
\def\jdcite{\cite}

\theoremstyle{plain}
\newtheorem{theorem}{Theorem}[section]
\newtheorem{lemma}[theorem]{Lemma}
\newtheorem{proposition}[theorem]{Proposition}

\theoremstyle{remark}

\title{Detecting complex multiplication}
\author{Jeffrey D. Achter}
\email{j.achter@colostate.edu}
\address{Department of Mathematics\\Colorado State University\\Fort
  Collins, CO 80523}

\begin{document}
\maketitle

\begin{abstract}
We give an efficient, deterministic algorithm to decide if two abelian varieties over a number field are isogenous.  From this, we derive an algorithm to compute the endomorphism ring of an elliptic curve over a number field.
\end{abstract}

In this paper, we answer two fundamental decision problems about
elliptic curves over number fields.  Specifically, we explain how to
detect whether two elliptic curves over a number field are isogenous,
and how to decide whether an elliptic curve has complex
multiplication.  These algorithms rely on Lemma
\ref{lemisog}, which actually applies to abelian varieties of any
dimension, and Proposition \ref{propcm}, respectively.

In each case, we answer a question about a variety over a number field
by examining its reduction at finitely many primes.   At this
level of generality, such a strategy is common in algorithmic number theory.
For example, a common method for computing modular polynomials --
that is, bivariate polynomials whose roots are $j$-invariants of
elliptic curves related by an isogeny of fixed degree -- is to perform the analogous computation
over various finite fields, and then to lift the result using the
Chinese remainder theorem.  In contrast, we will see that to answer
the decision problems posed here, one need not ever lift an object to
characteristic zero.

The engine driving the machines presented here is Faltings's paper on the Mordell
conjecture.  Milne observed in {\em Mathematical Reviews} that
Faltings ``seems to give an algorithm for deciding when two abelian
varieties over a number field are isogenous.''  In this paper, we
further refine the proof of \jdcite[Theorem 5]{faltings} to the point
where it literally yields an efficient algorithm for the isogeny
decision problem.

At a crucial stage in that argument, Faltings
shows that the isogeny class of $X$ is determined by the action of
$\gal(L/K)$ on $X[\ell](\bar K)$, where $[L:K]$ has effectively
bounded degree and ramification but is difficult to compute directly.
He therefore works with $\til L$, the compositum of all possible such
extensions of $K$, a large but still finite extension of $K$.  An
appeal to the Chebotarev density theorem guarantees that there is a
finite set of primes $T$ of $K$ such that $\st{ \artin{\til L/K}{\idp}
: \idp \in T} = \gal(\til L/K)$.  Therefore, $X$ and $Y$ are isogenous
if and only if the reductions $X_\idp$ and $Y_\idp$ are isogenous for
each $\idp\in T$.

We derive an algorithm for detecting isogeny by showing that it
suffices to use a set of primes $\idp$ with absolute norm smaller than
some constant  $B$.  Effective Chebotarev-type theorems
\jcite{bachsorenson,lagmontodly} let us calculate a
suitable $B$ solely in terms of the degree and ramification data of
$L$, without requiring recourse to the compositum $\til L$.  

Subsequently, we show how to use this result to test the hypothesis
that an elliptic curve $E$ has complex multiplication by a field $F$.
Briefly, after a finite extension of the base field, there exists an
elliptic curve $E'$ with complex multiplication by $F$.  Even without
computing $E'$ explicitly, we can use Lemma \ref{lemisog} to
detect whether $E$ and $E'$ are geometrically isogenous, and thus
check whether $E$ has complex multiplication by $F$.

%Our second main result is an effective, deterministic algorithm for
%deciding whether an elliptic curve over a number field has complex
%multiplication.  We give an effective upper bound $B$, depending on
%$E$ and $K$, such that $E$ has complex multiplication by a field $F$
%if and only if the same is true for all reductions $E_\idp$ where
%$\idp$ has norm at most $B$.

In the first section, we review literature concerning effective
Chebotarev density theorems, and explain how to use $\ell$-adic
representations to detect isogeny between abelian varieties.  The
reader may wish to skip Section \ref{subseccheb} on first reading, and
turn directly to Section \ref{subsecgalois}.

In the second section, we use these considerations to design
algorithms for elliptic curves over number fields.  In Section
\ref{subsecisog}, we describe an algorithm to determine whether two
elliptic curves are isogenous.  In Section \ref{subseccm}, we combine
the results of the previous section with new results on complex
multiplication to give an algorithm which decides whether a given
elliptic curve has complex multiplication.

Several improvements are available to improve the efficiency of these
methods. In the interest of streamlining the exposition, these
suggestions are gathered as a series of remarks in the final section.

I thank Siman Wong for helpful discussions.

\section{Background}

As discussed above, the method of this paper is to apply effective
Chebotarev bounds to Faltings's proof of the Tate conjecture in order
to construct efficient algorithms; we review these results in Sections
\ref{subseccheb} and \ref{subsecgalois}, respectively.

\subsection{Effective Chebotarev density theorems}
\label{subseccheb}

The Galois group of a finite extension of number fields $L/K$ is
generated by the Frobenius elements of primes of $L$ lying over primes
$\idp$ of $K$.  We collect here various results from the literature which place upper
bounds on the size of the primes necessary in order for their Artin
symbols to generate $\gal(L/K)$.    Throughout, we will use  $\grh$ to
highlight bounds which rely on the generalized Riemann hypothesis, and
$\nogrh$ to denote bounds which hold unconditionally.

For an extension of fields $L/K$, we let $\disc{L/K}$ denote the
discriminant and $\norm_{L/K}$ the norm map.  For a prime ideal $\idp$
of $K$, let $\kappa(\idp)$ be the residue field $\calo_K/\idp$ and let
$p_\idp$ be the characteristic of that field.

Let $S$ be a finite set of places of $K$ and $N$ a nonnegative
integer.   We will express our Chebotarev-type bounds in terms of the
following quantities:
\begin{align*}
\Delta^*(K,S,N) & \defeq
\abs{\disc{K/\rat}}^N ( N \cdot \prod_{\idp\in S}
p_\idp^{1-1/N})^{N\cdot [K:\rat]}\\
B_\lagod(K,S,N) & \defeq 70 \cdot (\log \Delta^*(K,S,N))^2 \\
B_\bach(K,S,N) & \defeq (4\log\Delta^*(K,S,N) + 2.5 N\cdot [K:\rat]+5)^2 \\
B_\grh(K,S,N) & \defeq \min\st{B_\lagod(K,S,N), B_\bach(K,S,N)}.
\end{align*}
Let $c_\nogrh$ be the effective constant $A_1$ of \jcite{lagmontodly}, and
let 
\begin{align*}
B_\nogrh(K,S,N) & \defeq 
\begin{cases}\Delta^*(K,S,N)^{c_\nogrh} & K\supsetneq \rat \\
2\Delta^*(K,S,N)^{c_\nogrh} & K = \rat
\end{cases}
.
\end{align*}
Finally, let
\begin{align*}
\calt_\grh(K,S,N) &\defeq \st{ \idp \subset K : \norm_{K/\rat} \idp \le
  B_\grh(K,S,N)\text{ and }\idp\not \in S} \\
\calt_\nogrh(K,S,N) &\defeq \st{ \idp \subset K : \norm_{K/\rat} \idp \le
  B_\bullet(K,S,N)\text{ and }\idp\not \in S}. \\
\end{align*}

%(Mnemonics: $B_\lagod$ for Lagarias and Odlyzko, $B_\bach$ for Bach
%and Sorenson, $B_\grh$ for generalized Riemann hypothesis, and
%$B_\nogrh$ for unconditional.)

%For $\bullet\in\st{\grh,\nogrh}$, let

\begin{lemma}\label{lemcheb} Let $K$ be a finite extension of $\rat$,
  and let $S\subset K$ be a finite set of prime ideals.  Let $L/K$ be
  a Galois extension with $[L:K]\le N$ unramified outside $S$.  For
  any $\sigma\in \gal(L/K)$, there exists $\idp \in
  \calt_\nogrh(K,S,N)$ and a prime $\idq$ of $L$ dividing $\idp$ such
  that $\frob_\idq = \sigma$.  If the generalized Riemann hypothesis
  holds, then $\idp$ may be taken in $\calt_\grh(K,S,N)$.
\end{lemma}

\begin{proof}
The statement combines several different effective Chebotarev
density theorems.  For a conjugacy class $C\subset\gal(L/K)$, each
gives an effective upper bound for the norm of the smallest prime
$\idp$ such that $\artin{L/K}{\idp} = C$, computed in terms of the
absolute discriminant of $L$.  By \jdcite[Proposition
5]{serrecheb}, $\abs{\disc{L/\rat}} \le \Delta^*(K,S,N)$; thus, in the
sequel, we may replace each occurrence of $\abs{\disc{L/\rat}}$ in
\jcite{bachsorenson,lagmontodly,odlylaga77} with $\Delta^*(K,S,N)$.

By \jdcite[Theorem 1.1]{lagmontodly}, any conjugacy class
$C\subset \gal(L/K)$ occurs as $\artin{L/K}{\idp}$ for some $\idp \in
\calt_\nogrh(K,S,N)$.  
Now suppose that the generalized Riemann hypothesis holds.  Lagarias
and Odlyzko prove \jcite{odlylaga77} that a bound of the form $B_\lagod$ suffices, and
Oesterle shows \jdcite[2.5]{serrecheb} that the constant is at most $70$.
The bound $B_\bach$ is obtained by Bach and Sorenson in
\jdcite[Theorem 5.1]{bachsorenson},  again under the assumption of the
generalized Riemann hypothesis.

Since the Frobenius elements $\frob_\idq$ of all primes lying over a
prime $\idp$ of $K$ form the conjugacy class $\artin{L/K}{\idp}$, the
result follows.
\end{proof}

\subsection{Abelian varieties and Galois modules}
\label{subsecgalois}

Let $X/K$ be an abelian variety, and let $\ell$ be a rational prime
such that $X$ has good reduction at all primes of $K$ lying over
$\ell$.  The $\ell$-adic Tate module of $X$ is $T_\ell(X)\defeq
\invlim n X[\ell^n](\bar K)$; let $V_\ell(X) \defeq
T_\ell(X)\tensor_\integ \rat$ be the rational Tate module.  Then
$T_\ell(X)$ is a $\integ_\ell$-representation of $\gal(\bar K/K)$,
while $V_\ell(X)$ is a $\rat_\ell$ representation of $\gal(\bar K/K)$.
It has long been known that these representations encode detailed
arithmetic information about $X$.

In fact, Faltings proves the Tate conjecture; the canonical map
$\End(X) \tensor_\integ \integ_\ell \ra
\End(T_\ell X)^{\gal(K)}$ is an isomorphism. Consequently
\jdcite[Corollary 2]{faltings} two abelian varieties are $X$ and $Y$ are
isogenous if and only if $V_\ell X$ and $V_\ell Y$ are isomorphic as
$\gal(\bar K/K)$-modules.

We denote the reduction of an abelian variety $X/K$ at a prime of good
reduction $\idp$
by $X_\idp$; it is an abelian variety over $\kappa(\idp)$.

The following result was proved by Serre \jdcite[8.3]{serrecheb} in the special case where
$\dim X = \dim Y = 1$ and $K = \rat$, but the the absolute constant
given there is ineffective.

\begin{lemma}\label{lemisog} Let $X$ and $Y$  be
  $g$-dimensional abelian varieties over a 
number field $K$.  Let $S$ be a set of places of $K$ containing all
primes of bad reduction of $X$ and $Y$, and let $\ell$ be a rational
prime which is relatively prime to each place of $S$.  Let
$\nu_g(\ell) = \abs{\gl_{2g}(\integ/\ell)}$.  Then $X$ and $Y$ are
isogenous if and only if $X_\idp$ and $Y_\idp$ are isogenous for all
$\idp\in \calt_\nogrh(K,S,\nu_g(\ell)^2)$.  If the generalized Riemann
hypothesis is true, then $X$ and $Y$ are isogenous if and only if
$X_\idp$ and $Y_\idp$ are isogenous for all $\idp\in
\calt_\grh(K,S,\nu_g(\ell)^2)$.
\end{lemma}

\begin{proof}
Our proof is closely modelled on that of \jdcite[Theorem 
  5]{faltings} and \jdcite[Theorem 23.7]{milneabvarnotes}.  Let $\calt =
\calt_\grh(K,S,\nu_g(\ell)^2)$ if the generalized Riemann hypothesis is
to be assumed, and let $\calt = \calt_\nogrh(K,S,\nu_g(\ell)^2)$
otherwise.    The key point is that the isogeny
class of an abelian variety over a number field is determined by the
Galois representation on its rational Tate module.  Lemma
  \ref{lemcheb} lets us detect
the isomorphism class of a Galois representation using only the
Frobenius elements over the finite set of primes $\calt$.

Let $\rho: \gal(\bar K/K) \ra \aut(T_\ell X) \cross \aut(T_\ell Y)$ be
the product representation.  Since $X$ and $Y$ both have good
reduction outside $S$, $\gal(\bar K/K)$ acts on $T_\ell X \cross
T_\ell Y$ via some quotient $\gal(E/K)$ with $E$ unramified outside
$S$.  Let $R$ be the subring of $\End(T_\ell X)\cross \End(T_\ell Y)$
generated over $\integ_\ell$ by $\st{ \rho(\sigma) : \sigma\in
\gal(E/K)}$.  We will show that $R$ is in fact
  generated, again over $\integ_\ell$, by the actions of $\frob_\idq$
  for primes $\idq$ of $E$ lying over $\idp\in \calt$.

By Nakayama's Lemma, it suffices to prove that these Frobenius
elements, acting on $(T_\ell X/\ell)\cross (T_\ell Y/\ell) =
X[\ell](\bar K)\cross Y[\ell](\bar K)$, generate $(R/\ell)\units$.  Now,
the action of $\gal(E/K)$ on $X[\ell](\bar K) \cross Y[\ell](\bar K)$
factors through 
$\gal(L/K)$, where $[L:K]$ is a finite Galois extension of degree at
most $\abs{ \aut(X[\ell](\bar K)) \cross \aut(Y[\ell](\bar K))} =
\nu_g(\ell)^2$.   By Lemma \ref{lemcheb}, $\st{\frob_\idq : \idq |
\idp\in\calt}= \gal(L/K)$.  Therefore,
$\st{\rho(\frob_\idq): \idq | \idp\in\calt }$ generates $R/\ell$ over
$\integ/\ell$, and this same set generates $R$ over $\integ_\ell$.

If $X_\idp$ and $Y_\idp$ are isogenous for some prime of good
reduction $\idp$, then \jdcite[Theorem 1]{tateendff} $V_\ell X$ and
$V_\ell Y$ are isomorphic as 
$\gal(\kappa(\idp))$-modules.  The hypothesis that $X_\idp$ and
$Y_\idp$ are isogenous for $\idp\in \calt$
implies that, for each $\frob_\idq$ with $\idq | \idp \in \calt$, $\tr( \frob_\idq | T_\ell X)
= \tr( \frob_\idq | T_\ell Y)$.  Extending $\integ_\ell$-linearly, we
have $\tr(\sigma | T_\ell X) = \tr(\sigma | T_\ell Y)$ for each
$\sigma\in \gal(\bar K/K)$, so that \jdcite[\S12.1, Proposition 3]{bourbakialgebre} $V_\ell X$ and $V_\ell Y$ are isomorphic
as $\gal(\bar K/K)$-modules.  By the Tate conjecture \jdcite[Corollary
2]{faltings} $X$ and $Y$ are isogenous.
\end{proof}

\section{Algorithms for elliptic curves}

\subsection{Detecting isogenous elliptic curves}
\label{subsecisog}

The isogeny class of an elliptic curve $E$ over a finite field
$\kappa$ is uniquely determined by $\abs{E(\kappa)}$.  Indeed, by
\jdcite[Theorem 1]{tateendff} the isogeny class of $E$ is determined by its
characteristic polynomial of Frobenius, which has the form $T^2 - aT +
\abs \kappa$.  Since the number of points on an elliptic curve with
such a characteristic polynomial is $\abs\kappa + 1 - a$, we see that
two elliptic curves over $\kappa$ are isogenous if and only if they
have the same number of points over $\kappa$.

Any efficient algorithm for counting points on elliptic curves over
finite fields, such as Schoof's method \jcite{schoofcount} which requires
$O(\log^9 \abs \kappa)$ bit operations, therefore yields an
efficient method for deciding if two elliptic curves are isogenous.

More generally, any efficient algorithm for computing the
action of Frobenius on $T_\ell X$ for a class of abelian varieties
$X$, such as Jacobians of hyperelliptic curves, can decide if two such
abelian varieties are isogenous.  (Note that, in dimension greater than
one, the action of Frobenius is {\em not} uniquely determined by its
trace.  Data such as the characteristic polynomial of Frobenius, rather than just
the trace of Frobenius, is required to detect the isogeny class of an
abelian variety over $\kappa$.)

We now turn our attention to number fields.   In
principle, Faltings's theorem affords us a choice of methods for
determining whether two elliptic curves $E_1$ and $E_2$ over a given
number field are
isogenous.  For instance, Masser and W\"ustholz
\jcite{masserwustholz90} use transcendence
theory to give an explicit  upper bound on the
minimal degree of an isogeny between two elliptic curves.  One could
then try to enumerate all curves related to $E_1$ by an isogeny of
given degree \jcite{velu}, and check if $E_2$ is isomorphic to any of
them.  One could also simply try to see if $E_1$ and $E_2$ satisfy a
modular equation of suitable degree.  Each of these operations carries
a nontrivial computational cost \jcite{agashelauter}.  Moreover, the best known {\em constant}
appearing in such degree bounds \jdcite[Th\'eor\`eme 1]{pellarin01} is
larger than $10^{61}$; such a method remains of theoretical, rather
than practical, interest.

Alternatively, an efficient algorithm follows from Lemma \ref{lemisog}.
Given two elliptic curves $E_1$ and $E_2$ over a common number field
$K$, compute the discriminant $\Delta_i$ of $E_i$, and thence the set
of primes $S_i$ for which $E_i$ has bad reduction.  If $S_1\not =
S_2$, then $E_1$ and $E_2$ are not isogenous \jdcite[Corollary 2]{serretate}.   Otherwise, choose a
rational prime $\ell$ relatively prime to each element of $S\defeq S_1
= S_2$. By Lemma \ref{lemisog}, $E_1$ and $E_2$ are
isogenous if and only if $E_{1,\idp}$ and $E_{2,\idp}$ are isogenous
for for each $\idp\in \calt(K, S,
(\ell^2-1)^2(\ell^2-\ell)^2)$; this last condition may be checked using
point-counting for each $E_{i,\idp}$.  (Again, if a method is
available for computing the characteristic polynomial of Frobenius,
then the same method works for detecting isogeny of abelian varieties
of dimension $g$; one simply computes at all primes with norm less
than $\nu_g(\ell)^2$.)

We remark that exhibiting infinitely many primes $\idp$ for which
$E_{1,\idp}$ and $E_{2,\idp}$ are isogenous 
does not prove that $E_1$ and $E_2$
are isogenous.  Indeed, suppose that $E_1$ and $E_2$ have complex
multiplication by distinct fields $F_1$ and $F_2$, respectively.  On
one hand, $E_1$ and $E_2$ are not isogenous, since the rational ring
of endomorphisms is an isogeny invariant.  On the other hand, we will
see below that 
$E_{1,\idp}$ and $E_{2,\idp}$ are both supersingular, and thus
isogenous, for all primes $\idp$ of $K$ for which $p_\idp$ is inert in
each extension $F_i/\rat$.

\subsection{Detecting complex multiplication}
\label{subseccm}

Let $E$ be an elliptic curve over a number field $K$.  The
endomorphism ring $\End(E)$ of $E$ is isomorphic either to $\integ$ or
to an order $\calo$ in a quadratic imaginary field, $F$.  In the
latter case, we say that $E$ has complex multiplication by $F$.  (More
generally, we will say that an elliptic curve over an arbitrary field
has complex multiplication by $F$ if its endomorphism ring contains an
order in $F$.)  

Elliptic curves with complex multiplication are prominent in primality
testing and cryptography \jcite{atkinmorain} and other aspects of
algorithmic number theory \jcite{cohencourse}.  Motivated
by this, one might seek an algorithm for determining whether a given
elliptic curve $E$ over a number field $K$ has complex multiplication.
In \jcite{charles}, the author describes two methods.  The first is a
probabilistic algorithm which runs in polynomial time in the inputs;
the second runs in deterministic polynomial time, but the constants
appearing in the analysis of the running time are ineffective. In this
section, we use Lemmas \ref{lemcheb} and \ref{lemisog} to give an
efficient, effective algorithm to determine whether an elliptic curve
has complex multiplication.  We start by collecting a body of facts
about elliptic curves with complex multiplication.  The subsequent
algorithm follows naturally from these observations.

Deuring investigated the relationship between the arithmetic of $F$
and the reductions $E_\idp$ at primes of $K$.  (For the
moment, we ignore primes of bad reduction.)  He proved (see
\jdcite[Exemple b]{tatehonda})
that $E_\idp$ is ordinary if and only if $p_\idp$, the rational prime
lying under $\idp$, splits in $F$.  Invoking the Chebotarev density
theorem for $F$, we see that $E$ has ordinary reduction at half the
primes of $K$, and supersingular reduction at the others.  

Conversely, if $\End(E)\iso \integ$, so that $E$ {\em does not} have
complex multiplication, then supersingular primes have density zero
\jdcite[IV-13, Exercise 1]{serreabelianladic}.
This basic observation leads to a probabilistic method, detailed in
\jcite{charles}, for checking whether an elliptic curve has complex
multiplication.  Broadly speaking, finding may primes of
supersingular reduction provides evidence for the hypothesis that $E$
has complex multiplication.

In the sequel, we will use Lemma \ref{lemisog} (and the accompanying
discussion at the end of Section \ref{subsecisog}) to describe a deterministic
algorithm to test whether an elliptic curve $E$ has complex
multiplication.  By this, we mean that the algorithm is guaranteed to
terminate after a finite, explicitly computable number of operations,
and that the output is a verifiable proof that $E$ does (or does not)
have complex multiplication.

It is convenient to assume that the elliptic curve of interest has
no automorphisms other than $\st{\pm 1}$, and that it has good
reduction everywhere.  The former condition is equivalent to the
assertion that $E$ {\em does not} have complex multiplication by
$\rat(\sqrt{-1})$ or $\rat(\sqrt{-3})$, which is easily verified by
checking that $j(E)\not\in\st{0,1728}$.

The latter condition holds, possibly after a finite extension of the
base field, for an elliptic curve with complex multiplication.  (This
assertion is equivalent to the result of Weber
\jdcite[C.11.2.a]{silvermanaec} that an elliptic curve with complex
multiplication has integral $j$-invariant.)  Concretely, let $\caln$
be the product of all primes of bad reduction of $E$, and let $K_1 =
K(\sqrt{\caln})$.  Suppose that $E$ has complex multiplication by a
field whose only roots of unity are $1$ and $-1$.  A special case of
\jdcite[Theorem 7]{serretate} shows that $E_{K_1}$ has good reduction
at all places of $K_1$.

Henceforth, we will assume that $E/K$ has everywhere good reduction.

We now show that it is easy to find a prime of ordinary reduction for
$E$, and thereby find a candidate ring of endomorphisms for $E$.

\begin{lemma}\label{lemord}  Let $E/K$ be an elliptic curve over a number field with
complex multiplication and good reduction everywhere.  Let $S$ be the set of rational primes
ramified in the extension $S/\rat$.  Then there exists a prime $\idp$ of $K$ lying over
a rational prime $p$ with $p \le B_{\nogrh}(\rat, S, 2)$ such
that $E$ has good, ordinary reduction at $\idp$.  If the generalized
Riemann hypothesis is true, then $p$ may be taken less than or equal
to $B_\grh(\rat, S_\rat, 2)$.
\end{lemma}
\begin{proof}
If $\End(E)\tensor\rat$ is isomorphic to a quadratic imaginary field
$F$, then $K$ necessarily contains $F$ \jdcite[Theorem 3.1.1]{langcm}.  In particular, the support of
the discriminant of $F$ over $\rat$ is contained in the support of the
discriminant of $K$ over $\rat$, so that $F$ is unramified outside
$S$.  Moreover, $E$ has ordinary reduction at a prime $\idp$ over the
rational prime $p$ if and only if $p$ splits in $F$. The Chebotarev density theorem
(Lemma \ref{lemcheb}) guarantees the existence of such a $p$ with $p
\le B_{\nogrh}(\rat, S_\rat, 2)$.  (If the generalized Riemann
hypothesis is true, then $p$ may be taken to be at most
$B_{\grh}(\rat, S_\rat, 2)$.)
\end{proof}

Since supersingular primes have density zero for an elliptic curve {\em
without} complex multiplication, it seems unlikely that one would
encounter an $E/K$ without a small (in the sense of Lemma
\ref{lemord}) ordinary prime.  Still, if this were to happen, one
could then conclude that the elliptic curve had endomorphism ring
equal to $\integ$.

Let $\idp$ be a prime of ordinary reduction of $E$.  Then the ring
$\integ[\frob_\idp] \subseteq \End(E_\idp)$ is isomorphic to an order
in a quadratic imaginary field $F$.  (To see this, use the result of
Deuring \jdcite[Theorem V.3.1]{silvermanaec}, paralleling the result in
characteristic zero, that the 
endomorphism ring of an ordinary elliptic curve is either $\integ$ or
an order in a quadratic imaginary field.  Moreover, the Frobenius
endomorphism cannot have a real conjugate \jdcite[Exemple
a]{tatehonda}, so it must actually generate a quadratic imaginary field.)
Moreover, given the number of
points on $E_\idp$, one can determine the field $F = \Frac
\integ[\frob_\idp]$.  This is a candidate field of (rational)
endomorphisms of $E$, and we show in Proposition \ref{mainprop} how to
test the hypothesis that $E$ truly does have complex multiplication by
$F$.   (At this stage of the calculation, one knows that $E$ has
complex multiplication by {\em some}
quadratic imaginary field if and only if it has complex multiplication
by $F$.)

It is known \jdcite[Corollary C.11.1]{silvermanaec} that there is a finite extension $K'$ of $K$ and an
elliptic curve $E'/K'$ with complex multiplication by $F$.  At this point,
one could simply compute the $j$-invariant of $E'$ and check whether
$j(E)$ and $j(E')$ are conjugate under $\gal(\rat)$.  However, to
compute the polynomial over $\rat$ which $j(E')$ satisfies takes time
$O(\abs{\disc{F/\rat}}^2 (\log\abs{\disc{F/\rat}})^2)$
\jcite{atkinmorain} (see also \jdcite[7.6]{cohencourse}).

Now, if one could construct $E'$ efficiently, one could use Lemma 
\ref{lemisog} to test whether $E$ and $E'$ are isogenous, since 
isogenous elliptic curves have commensurable rings of endomorphisms. 
Even without knowing $E'$ explicitly, however, we can efficiently test
whether the two curves are (geometrically) isogenous.

\begin{proposition}\label{mainprop}\label{propcm} Suppose that $E/K$ has good reduction everywhere.  Let
$F$ be a quadratic imaginary subfield of $K$ whose only roots of unity are $-1$
and $1$, and let $h^*(F) = 2\sqrt{\disc{F/\rat}}/\pi$. 
Then
  $E$ has complex multiplication by $F$ if and only if for each prime
  $\idp\in K$ lying over a rational prime $p$ with $\norm \idp \le
B_\nogrh(K,\emptyset,h^*(F)\nu_2(2)^2)$, either:
\begin{itemize} 
\item $E_\idp$ is supersingular and $F$ is inert or ramified at $p$, or
\item $\End(E_\idp)\tensor\rat\iso F$, and $F$ is split at $p$.
\end{itemize}
If the generalized Riemann hypothesis is true, then it suffices to 
consider those primes with norm at most $B_\grh(K,\emptyset,h^*(F)\nu_2(2)^2)$. 
\end{proposition}

\begin{proof}
If the generalized Riemann hypothesis is to be assumed, we write
$B$ for $B_\grh$ and $\calt$ for $\calt_\grh$; otherwise, these
symbols denote $B_\nogrh$ and $\calt_\nogrh$, respectively.  Note that
the statement is equivalent to the assertion that $E$ has complex
multiplication by $F$ if and only if the same is true of $E_\idp$ for
each prime $\idp \in B(K,\emptyset, h^*(F)\nu_2(2)^2)$.  Since there
is a natural inclusion $\End(E)\inject \End(E_\idp)$ for each prime
$\idp$ \jdcite[Theorem 2.3.2]{langcm}, if $E$ has complex multiplication
by $F$ then the same is true of $E_\idp$ for each prime $\idp$, and in
particular for those in $B(K,\emptyset, h^*(F)\nu_2(2)^2)$.

Having secured this, we focus on the converse.  There exists an
elliptic curve over a field $K'$ with complex multiplication by $F$ if
and only if $K'$ contains the Hilbert class field of $F$
\jdcite[Theorem C.11.2]{silvermanaec}.  Moreover,
since the only roots of unity in $F$ are $\st{\pm 1}$, we may assume
that $E'$ has good reduction everywhere \jdcite[Theorem 9]{serretate}.

Therefore, let $K'$ be the compositum of $K$ and the Hilbert class
field of $F$, and let $E'/K'$ be an elliptic curve with everywhere
good reduction and complex multiplication by $F$.  The original
elliptic curve $E$ has complex multiplication by $F$ if and only if
$E_{K'}$ and $E'$ are isogenous over some finite extension of $K'$.
(An analytic construction, as in \jdcite[C.11]{silvermanaec}, shows
that $E$ and $E'$ are isogenous over ${\mathbb C}$; this isogeny must then
descend to some finite $K'/K$ \jdcite[Theorem II.2.2]{silvermanadvtop}.)
Equivalently, $E$ has complex multiplication by $F$ if and only if
$E_{K'}$ is isogenous to some twist of $E'$.

Let $N = \nu_2(2)^2$ and suppose that,  for all primes $\idq\in \calt(K', \emptyset, N)$,
$E_\idq$ has complex multiplication by $F$.  Then $E_\idq$ and
$E'_\idq$ are isogenous up to a quadratic 
twist, and there exists a twist $E''$ of $E'$ such that
$E_\idq$ and $E''_\idq$ are isogenous for all $\idq\in
\calt(K',\emptyset,N)$.  By Lemma \ref{lemisog}, $E$ and $E''$ are
isogenous, and thus $E$ has complex multiplication by $F$.

If the prime $\idq$ of $K'$ lies over the prime $\idp$ of $K$, then
$\norm_{K'/\rat} (\idq) \le \norm_{K/\rat}(\idp)$.  In particular,
each prime $\idq\in \calt(K',\emptyset, N)$ lies over a prime $\idp$
with $\norm_{K/\rat}\idp \le B(K',\emptyset, N)$.  Moreover, $E_\idq$
is the base change $E_\idp \cross \kappa(\idq)$, and thus $E_\idq$ has complex
multiplication by $F$ if and only if $E_\idp$ does.  We have thus
shown that $E$ has complex multiplication by $F$ if and only if the
same is true for each reduction $E_\idp$ with $\norm_{K/\rat}\idp \le
B(K',\emptyset, N)$.  Now, $K'$ is an unramified extension of $K$ of
degree at most $h^*(F)$.  Therefore,
\begin{align*}
\Delta^*(K',\emptyset, N) &= \abs{\Delta_{K'/\rat}}^N
N^{N\cdot[K':\rat]} \\ &= \abs{\Delta_{K/\rat}^{[K':K]}}^N \cdot
N^{N\cdot[K':K]\cdot[K:\rat]} \\ &= \abs{\Delta_{K/\rat}}^{N\cdot
[K':K]}(N^{N\cdot[K':K]})^{[K:\rat]} \\ & \le \Delta^*(K,\emptyset,
h^*(F)N),
\end{align*}
and the result follows.
\end{proof}

Taken together, the results of this subsection suggest the following
algorithm for determining whether an elliptic curve $E$ has complex
multiplication.  First, check whether $j(E)\in \st{0,1728}$; if so,
the answer is yes; if not, one continues.  Second, construct
$K(\sqrt{\caln})$, and verify that $E/K(\sqrt{\caln})$ has good
reduction everywhere; if this fails, then \jdcite[Theorem 7]{serretate} $E$
does not have complex multiplication.  Otherwise, replace $K$ with
$K(\sqrt{\caln})$, and use Lemma \ref{lemord} to find a candidate
field $F$ of endomorphisms.  Finally, Proposition \ref{mainprop}
allows us to test efficiently if $E/K$ has complex multiplication by
$F$.

\paragraph{\null}

\section{Algorithmic considerations}
We close with some remarks which may allow more efficient
implementation of these algorithms.

\subsection{Chebotarev density theorem}
It is sometimes possible to improve the bounds given in Lemma \ref{lemcheb}.
In special cases where $[K:\rat]$ and $N$ are both small, \jdcite[Table
1]{bachsorenson} provides even tighter bounds for the norm of the
smallest prime ideal with given Artin symbol.  

Moreover, if either the bound $B_\bach$ or $B_\nogrh$  is used, then it suffices to consider those primes
$\idp$ of $K$ with norm a rational prime.

\subsection{Bounds in Lemma \ref{lemisog}}
The term $\nu(g) = \abs{\gl_{2g}(\integ/\ell)}$ arises in the proof of
Lemma \ref{lemisog} as the size of the automorphism group of
$X[\ell]$.  If $X$ further comes equipped with a polarization over $K$
of degree prime to $\ell$, then the action of $\gal(\bar K/K)$ on
$T_\ell X$ commutes with the induced symplectic pairing.  Therefore,
if one further makes the assumption in Lemma \ref{lemisog} that $X$
and $Y$ admit polarizations over $K$ of degree relatively prime to
$\ell$, then $\nu_g(\ell)$ may be replaced by
$\abs{\gsp_{2g}(\integ/\ell)}$.

\subsection{Candidate fields of complex multiplication}
In Section \ref{subseccm}, we suggested using Lemma \ref{lemord} to
find a candidate ring of endomorphisms of $E$.  Alternatively, if
$\End_K(E)$ is an order in a quadratic imaginary field then $\End_K(E)
\subseteq K$, so that $\End_K(E)\tensor\rat$ is a quadratic imaginary
subfield of $K$.  Therefore, one can enumerate each such subfield
$F_i$ of $K$, and apply Proposition \ref{mainprop} to each; $E$ has
complex multiplication by some field if and only if it has complex
multiplication by one of the $F_i$.

\subsection{Bounds in Proposition \ref{mainprop}}
The proof of Proposition \ref{mainprop} shows that it suffices to consider
those primes of $K$ with norm at most, e.g.,
\[
70\cdot( h^*(F) \log \Delta^*(K,\emptyset, \nu_2(2)^2))^2
\]
if the Lagarias and Odlyzko bound is to be used; the analogous
improvement may be made in each of the other bounds, as well.
Moreover, one can replace the (perhaps pessimistic) bound $h^*(F)$
with the actual class number of $F$; this class number can be computed
in time $O(\abs{\disc{F/\rat}}^{1/4+\epsilon})$, and even
$O(\abs{\disc{F/\rat}}^{1/5+\epsilon})$ if the generalized Riemann
hypothesis is assumed \jdcite[5.4]{cohencourse}.  Finally, if one could
verify that the Hilbert class field of $F$ is already contained in
$K$, then one would know (in the notation of the proof) that $K= K'$,
and one could replace $h^*(F)$ with $1$.  However, it is not clear to
the author how to verify this condition, short of actually computing
the Hilbert class field.

\subsection{Detecting potential complex multiplication}
It is not hard to adopt the observations of Section \ref{subseccm} to
test whether $E$ {\em potentially} has complex multiplication, in the
sense that $\End_{\bar K}(E)\tensor\rat$ is a quadratic imaginary
field.  One needs to replace Lemma \ref{lemord} with an upper bound
for the size of a prime of ordinary reduction; as noted there, we
expect in practice that it is quite easy to find such a prime.  This
generates a candidate field $F$ of rational endomorphisms for $E_{\bar
K}$.  One can then apply Proposition \ref{mainprop} to check whether
$E_{K.F}$ has complex multiplication by $F$.

\bibliographystyle{hplain} 
\bibliography{ecrefs}

\end{document}